\documentclass[12pt]{article}
\usepackage{amssymb,amsmath}
\usepackage{hyperref}
\usepackage{graphicx}
\usepackage{color}
\usepackage{chngcntr}

\begin{document}

\title{\LARGE\bf A rational approximation of \\
the Dawson\text{'}s integral for efficient computation of the complex \\
error function}

\author{
\normalsize\bf S. M. Abrarov\footnote{\scriptsize{Dept. Earth and Space Science and Engineering, York University, Toronto, Canada, M3J 1P3.}}
\, and B. M. Quine$^{*}$\footnote{\scriptsize{Dept. Physics and Astronomy, York University, Toronto, Canada, M3J 1P3.}}}

\date{November 23, 2017}
\maketitle

\begin{abstract}
In this work we show a rational approximation of the Dawson\text{'s} integral that can be implemented for high accuracy computation of the complex error function in a rapid algorithm. Specifically, this approach provides accuracy exceeding $\sim {10^{ - 14}}$ in the domain of practical importance $0 \le y < 0.1 \cap \left| {x + iy} \right| \le 8$. A Matlab code for computation of the complex error function with entire coverage of the complex plane is presented.
\vspace{0.25cm}
\vspace{0.25cm}
\end{abstract}

\section{Introduction}

The complex error function also widely known as the Faddeeva function can be defined as \cite{Faddeyeva1961, Gautschi1970, Abramowitz1972, Armstrong1972, Schreier1992, Weideman1994}
\begin{equation}\label{eq_1}
w\left( z \right) = {e^{ - {z^2}}}\left( {1 + \frac{{2i}}{{\sqrt \pi  }}\int\limits_0^z {{e^{{t^2}}}dt} } \right),
\end{equation}
where $z = x + iy$ is the complex argument. It is a solution of the following differential equation \cite{Schreier1992}
\[
w'\left( z \right) + 2zw\left( z \right) = \frac{{2i}}{{\sqrt \pi  }},
\]
where the initial condition is given by $w\left( 0 \right) = 1.$

The complex error function is a principal in a family of special functions. The main functions from this family are the Dawson\text{'}s integral, the complex probability function, the error function, the Fresnel integral and the normal distribution function.

The Dawson\text{'}s integral is defined as \cite{Cody1970, McCabe1974, Rybicki1989, Boyd2008, Abrarov2015a, Nijimbere2017}
\begin{equation}\label{eq_2}
{\rm{daw}}\left( z \right) = {e^{ - {z^2}}}\int\limits_0^z {{e^{{t^2}}}dt}.
\end{equation}
It is not difficult to obtain a relation between the complex error function and the Dawson\text{'s} integral. In particular, comparing right sides of equations \eqref{eq_1} and \eqref{eq_2} immediately yields
\begin{equation}\label{eq_3}
w\left( z \right) = {e^{ - {z^2}}} + \frac{{2i}}{{\sqrt \pi  }}{\rm{daw}}\left( z \right).
\end{equation}

Another closely related function is the complex probability function. In order to emphasize the continuity of the complex probability function at $\forall y \in \mathbb{R} $, it may be convenient to define it in form of principal value integral \cite{Armstrong1972, Schreier1992, Weideman1994}
\begin{equation}\label{eq_4}
W\left( z \right) = PV\frac{i}{\pi }\int\limits_{ - \infty }^\infty  {\frac{{{e^{ - {t^2}}}}}{{z - t}}dt}
\end{equation}
or
\[
W\left( {x,y} \right) = PV\frac{i}{\pi }\int\limits_{ - \infty }^\infty  {\frac{{{e^{ - {t^2}}}}}{{\left( {x + iy} \right) - t}}dt}.
\]
The complex probability function has no discontinuity at $y = 0$ and $x = t$ since according to the principal value we can write
\begin{equation}\label{eq_5}
\lim W\left( {x,y \to 0} \right) = {e^{ - {x^2}}} + \frac{{2i}}{{\sqrt \pi  }}{\rm{daw}}\left( x \right),
\end{equation}
where $x = {\mathop{\rm Re}\nolimits} \left[ z \right]$.

There is a direct relationship between the complex error function \eqref{eq_1} and the complex probability function \eqref{eq_4}. In particular, it can be shown that these functions are actually same on the upper half of the complex plain \cite{Armstrong1972, Schreier1992}
\begin{equation}\label{eq_6}
W\left( z \right) = w\left( z \right), \qquad {\mathop{\rm Im}\nolimits} \left[ z \right] \ge 0.
\end{equation}
Separating the real and imaginary parts of the complex probability function \eqref{eq_4} leads to
$$
K\left( {x,y} \right) = PV\frac{y}{\pi }\int\limits_{ - \infty }^\infty  {\frac{{{e^{ - {t^2}}}}}{{{y^2} + {{\left( {x - t} \right)}^2}}}dt,}
$$
and
\begin{equation}\label{eq_7}
L\left( {x,y} \right) = PV\frac{1}{\pi }\int\limits_{ - \infty }^\infty  {\frac{{{e^{ - {t^2}}}\left( {x - t} \right)}}{{{y^2} + {{\left( {x - t} \right)}^2}}}dt,}
\end{equation}
respectively, where the principal value notations emphasize that these functions have no discontinuity at $y = 0$ and $x = t$. In particular, in accordance with equation \eqref{eq_3} we can write
$$
K\left( {x,y = 0} \right) \equiv \lim K\left( {x,y \to 0} \right) = {e^{ - {x^2}}}
$$
and
$$
L\left( {x,y = 0} \right) \equiv \lim L\left( {x,y \to 0} \right) = \frac{{2}}{{\sqrt \pi  }}\,{\rm{daw}}\left( x \right).
$$
As it follows from the identity \eqref{eq_6}, for non-negative $y$ we have
\begin{equation}\label{eq_8}
w\left( {x,y} \right) = K\left( {x,y} \right) + iL\left( {x,y} \right), \qquad y \ge 0.
\end{equation}

The real part $K\left( {x,y} \right)$ of the complex probability function is commonly known as the Voigt function that is widely used in many disciplines of Applied Mathematics \cite{Srivastava1987, Srivastava1992, Pagnini2010}, Physics \cite{Armstrong1972, Armstrong1967, Letchworth2007, Edwards1992, Quine2002, Christensen2012, Berk2013, Quine2013}, Astronomy \cite{Emerson1996} and Information Technology \cite{Farsad2015}. Mathematically, the Voigt function $K\left( {x,y} \right)$ represents a convolution integral of the Gaussian and Cauchy distributions \cite{Schreier1992, Armstrong1967, Letchworth2007}. The Voigt function is widely used in spectroscopy as it describes quite accurately the line broadening effects \cite{Armstrong1972, Edwards1992, Quine2002, Christensen2012, Berk2013, Quine2013, Prince2004}.

Although the imaginary part $L\left( {x,y} \right)$ of the complex probability function also finds many practical applications (see for example \cite{Balazs1969, Chan1986}), it has no a specific name. Therefore, further we will regard this function simply as the $L$-function.

Other associated functions are the error function of complex argument \cite{Abramowitz1972, Schreier1992}
$$
{\rm{erf}}\left( z \right) = 1 - {e^{ - {z^2}}}w\left( {iz} \right),
$$
the plasma dispersion function \cite{Fried1961}
$$
Z\left( z \right) = PV\frac{1}{{\sqrt \pi  }}\int\limits_{ - \infty }^\infty  {\frac{{{e^{ - {t^2}}}}}{{t - z}}dt = i\sqrt \pi  w\left( z \right)}
$$
the Fresnel integral \cite{Abramowitz1972}
$$
\begin{aligned}
F\left( z \right) &= \int\limits_{ 0 }^z  {{e^{i\left( {\pi /2} \right){t^2}}}dt} \\
 &= \left( {1 + i} \right)\left[ {1 - {e^{i\left( {\pi /2} \right){z^2}}}w\left( {\sqrt \pi  \left( {1 + i} \right)z/2} \right)} \right]/2
\end{aligned}
$$
and the normal distribution function \cite{Weisstein2003}
$$
\begin{aligned}
\Phi \left( z \right) &= \frac{1}{{\sqrt {2\pi } }}\int\limits_0^z {{e^{ - {t^2}/2}}dt = \frac{1}{2}{\rm{erf}}\left( {\frac{z}{{\sqrt 2 }}} \right)} \\
 &= \frac{1}{2}\left[ {1 - {e^{ - {z^2}/2}}w\left( {\frac{{iz}}{{\sqrt 2 }}} \right)} \right].
\end{aligned}
$$

Due to a remarkable identity of the complex error function \cite{McKenna1984, Zaghloul2011}
\[
w\left( { - z} \right) = 2{e^{ - {z^2}}} - w\left( z \right),
\]
it is sufficient to consider only I and II quadrants in order to cover the entire complex plane (see Appendix A). This property can be explicitly observed by rearranging this identity in form
$$
w\left( { \pm x, - \left| y \right|} \right) = 2{e^{ - {{\left( { \mp x + i\left| y \right|} \right)}^2}}} - w\left( { \mp x, + \left| y \right|} \right).
$$
Thus, due to this remarkable identity of the complex error function we can always avoid computations involving negative argument $y$ by using simply the right side of this equation. Therefore, we will always assume further that $y \ge 0$.

Since the integral on right side of the equation \eqref{eq_1} has no analytical solution, the computation of the complex error function must be performed numerically. There are many approximations of the complex error function have been reported in the scientific literature \cite{Weideman1994, Zaghloul2011, Press1992, Chiarella1968, Abrarov2014, Abrarov2015b, Abrarov2011}. However, approximations based on rational functions may be more advantageous to gain computational efficiency due to no requirement for nested loop involving decelerating exponential or trigonometric functions dependent upon input the parameters $x$ and $y$. Although the rational approximations of the complex error function are rapid and provide high-accuracy over the most area of the complex plain, all of them have the same drawback; their accuracy deteriorates with decreasing parameter $y$ (in fact, the computation of the complex error function at small $y <  < 1$ is commonly considered problematic for high-accuracy computation of the Voigt/complex error function not only in the rational approximations \cite{Armstrong1967, Amamou2013, Abrarov2015c}). In this work we propose a new rational approximation of the Dawson\text{'s} integral and show how its algorithmic implementation effectively resolves such a problem in computation of the complex error function at small $y <  < 1$.

\section{Derivation}

In our recent publication we have shown that a new methodology of the Fourier transform results in a rational approximation \cite{Abrarov2015d}
\begin{equation}\label{eq_9}
w\left( z \right) \approx i\frac{{2{h}{e^{{\sigma ^2}}}}}{z +i\sigma } + \sum\limits_{n = 1}^N {\frac{{{A_n} - i\left( {z + i\sigma } \right){B_n}}}{{C_n^2 - {{\left( {z + i\sigma } \right)}^2}}}},
\end{equation}
where the corresponding expansion coefficients are given by
$$
{A_n} = 8\pi h^2n{e^{{\sigma ^2} - {{\left( {2\pi {h}n} \right)}^2}}}\sin \left( {4\pi {h}n\sigma } \right),
$$
$$
{B_n} = 4{h}{e^{{\sigma ^2} - {{\left( {2\pi {h}n} \right)}^2}}}\cos \left( {4\pi {h}n\sigma } \right),
$$
$$
{C_n} = 2\pi {h}n
$$
and $h$ is the step. A brief description for derivation of the equation \eqref{eq_9} is presented in Appendix B. In algorithmic implementation it is convenient to rewrite the equation above in form
\begin{equation}\label{eq_10}
\begin{aligned}
\psi \left( z \right) &= i\frac{{2{h}{e^{{\sigma ^2}}}}}{z} + \sum\limits_{n = 1}^N {\frac{{{A_n} - iz{B_n}}}{{C_n^2 - {z^2}}}} \\
 &\Rightarrow w\left( z \right) \approx \psi \left( {z + i\sigma } \right).
\end{aligned}
\end{equation}

Using the residue calculus, from approximation \eqref{eq_9} we can obtain the rational approximation for the $L$-function (see Appendix C)
\begin{equation}\label{eq_11}
L\left( x \right) \approx x \left[ \frac{{2{e^{{\sigma ^2}}}h}}{{{x^2} + {\sigma ^2}}} + \sum\limits_{n = 1}^N {\frac{{2{A_n}\sigma  + {B_n}\left( {{x^2} + {\sigma ^2} - C_n^2} \right)}}{{C_n^4 + 2C_n^2\left( {{\sigma ^2} - {x^2}} \right) + {{\left( {{x^2} + {\sigma ^2}} \right)}^2}}}} \right].
\end{equation}
As $L\left( {x,y = 0} \right) = {\mathop{\rm Im}\nolimits} \left[ {w\left( {x,y = 0} \right)} \right]$ from the identity \eqref{eq_3} it follows that
\[
{\rm{daw}}\left( x \right) = \frac{{\sqrt \pi  }}{2}{\mathop{\rm Im}\nolimits} \left[ {w\left( x \right)} \right] = \frac{{\sqrt \pi  }}{2}L\left( x \right).
\]
Consequently, according to this identity and equation \eqref{eq_11} we can write
\begin{equation}\label{eq_12}
{\rm{daw}}\left( x \right) \approx \frac{{\sqrt \pi  x}}{2}\left[ {\frac{{2{e^{{\sigma ^2}}}h}}{{{x^2} + {\sigma ^2}}} + \sum\limits_{n = 1}^N {\frac{{2{A_n}\sigma  + {B_n}\left( {{x^2} + {\sigma ^2} - C_n^2} \right)}}{{C_n^4 + 2C_n^2\left( {{\sigma ^2} - {x^2}} \right) + {{\left( {{x^2} + {\sigma ^2}} \right)}^2}}}} } \right].
\end{equation}
Since the argument $x$ is real, it would be reasonable to assume that this equation remains valid for a complex argument $z = x + iy$ as well if the condition $y <  < 1$ is satisfied. Therefore, substituting the rational approximation of the Dawson\text{'s} integral \eqref{eq_12} into identity \eqref{eq_3} results in
\small
\begin{equation}\label{eq_13}
\begin{aligned}
w\left( z \right) \approx {e^{ - {z^2}}} &+ iz\left[ \frac{{2{e^{{\sigma ^2}}}h}}{{{z^2} + {\sigma ^2}}} \right. \\
&\left. + \sum\limits_{n = 1}^N {\frac{{2{A_n}\sigma  + {B_n}\left( {{z^2} + {\sigma ^2} - C_n^2} \right)}}{{C_n^4 + 2C_n^2\left( {{\sigma ^2} - {z^2}} \right) + {{\left( {{z^2} + {\sigma ^2}} \right)}^2}}}}  \right], \qquad	{\mathop{\rm Im}\nolimits} \left[ z \right] <  < 1.
\end{aligned}
\end{equation}
\normalsize

The representation of the approximation \eqref{eq_13} can be significantly simplified as given by
\begin{equation}\label{eq_14}
\begin{aligned}
w\left( z \right) &\approx {e^{ - {z^2}}} + 2ih{e^{{\sigma ^2}}}z\theta \left( {{z^2} + {\sigma ^2}} \right)\\
 &\Rightarrow \theta \left( z \right) \buildrel \Delta \over = \frac{1}{z} + \sum\limits_{n = 1}^N {\frac{{{\alpha _n} + {\beta _n}\left( {z - {\gamma _n}} \right)}}{{4{\sigma ^2}{\gamma _n} + {{\left( {{\gamma _n} - z} \right)}^2}}}}, \quad  {\mathop{\rm Im}\nolimits} \left[ z \right] <  < 1,
\end{aligned}
\end{equation}
\normalsize
where the expansion coefficients are
$$
{\alpha _n} = \frac{\sigma }{{h{e^{{\sigma ^2}}}}}{A_n} = 8\pi h n\sigma {e^{ - {{\left( {2\pi hn} \right)}^2}}}\sin \left( {4\pi hn\sigma } \right),
$$
$$
{\beta _n} = \frac{1}{{2h{e^{{\sigma ^2}}}}}{B_n} = 2{e^{ - {{\left( {2\pi hn} \right)}^2}}}\cos \left( {4\pi hn\sigma } \right)
$$
and
$$
{\gamma _n} = C_n^2 = {\left( {2\pi hn} \right)^2}.
$$

As we can see from equation \eqref{eq_14} only the $\theta $-function needs a nested loop in multiple summation. Therefore, most time required for computation of the complex error function is taken for determination of the $\theta $-function. However, this approximation is rapid due to its simple rational function representation. Although the first term ${e^{ - {z^2}}}$ in approximation \eqref{eq_14} is an exponential function dependent upon the input parameter $z$, it does not decelerate the computation since this term is not nested and calculated only once. Consequently, the approximation \eqref{eq_14} is almost as fast as the rational approximation \eqref{eq_10}.

\section{Implementation and error analysis}

The relatively large area at the origin of complex plane is the most difficult for computation with high-accuracy. According to Karbach {\it{et al.}} \cite{Karbach2014} the newest version of the {\it{RooFit}} package, written in C/C++ code, utilizes the equation (14) from the work \cite{Abrarov2011} in order to cover accurately this area, shaped as a square with side lengths equal to $2 \times 8 = 16$.

In the algorithm we have developed, instead of the square we apply a circle with radius $8$ centered at the origin. This circle separates the complex plane into the inner and outer parts. Only three approximations can be applied to cover the entire complex plain with high-accuracy. The outer part of the circle is an external domain while the inner part of the circle is an internal domain consisting of the primary and secondary subdomains. These domains are schematically shown in Fig. 1.

External domain is determined by boundary $\left| {x + iy} \right| > 8$ and represents the outer area of the circle as shown in Fig. 1. In order to cover this domain we apply the truncated Laplace continued fraction \cite{Gautschi1970, Jones1988, Poppe1990a, Poppe1990b}
$$
w\left( z \right) = \frac{{{\mu _0}}}{{z - }}\frac{{1/2}}{{z - }}\frac{1}{{z - }}\frac{{3/2}}{{z - }}\frac{2}{{z - }}\frac{{5/2}}{{z - }}\frac{3}{{z - }}\frac{{7/2}}{{z - }}...,	\qquad {\mu _0} = i/\sqrt \pi .
$$
that provides a rapid computation with high-accuracy. It should be noted, however, that the accuracy of the Laplace continued fraction deteriorates as $\left| z \right|$ decreases.

The inner part of the circle bounded by $\left| {x + iy} \right| \le 8$ is divided into two subdomains. Most area inside the circle is occupied by primary subdomain bounded by $y \ge 0.1 \cap \left| {x + iy} \right| \le 8$. For this domain we apply the rational approximation \eqref{eq_10} that approaches the limit of double precision when we take $N = 23$, $\sigma  = 1.5$ and $h = 6/\left( {2\pi N} \right)$ (see our recent publications \cite{Abrarov2015d} and \cite{Abrarov2015e} for description in determination of the parameter $h$).

The secondary subdomain within the circle is bounded by narrow band $0 \le y < 0.1 \cap \left| {x + iy} \right| \le 8$ (see Fig. 1 along $x$-axis). The rational approximation \eqref{eq_10} sustains high-accuracy within all domain $0 \le x \le 40,000 \cap {10^{ - 4}} \le y \le {10^2}$ required for applications using the HITRAN molecular spectroscopic database \cite{Rothman2013}. However, at $y < {10^{ - 6}}$ its accuracy deteriorates by roughly one order of the magnitude as $y$ decreases by factor of $10$. The proposed approximation \eqref{eq_14} perfectly covers the range $y <  < 1$. In particular, at $N = 23$, $\sigma  = 1.5$ and $h = 6/\left( {2\pi N} \right)$ this approximation also approaches the limit of double precision as the parameter $y$ tends to zero.

In order to quantify the accuracy of the algorithm we can use the relative errors defined for the real and imaginary parts as follows
$$
{\Delta _{{\mathop{\rm Re}\nolimits} }} = \left| {\frac{{{\mathop{\rm Re}\nolimits} \left[ {{w_{ref}}\left( {x,y} \right)} \right] - {\mathop{\rm Re}\nolimits} \left[ {w\left( {x,y} \right)} \right]}}{{{\mathop{\rm Re}\nolimits} \left[ {{w_{ref}}\left( {x,y} \right)} \right]}}} \right|
$$
and
$$
{\Delta _{{\mathop{\rm Im}\nolimits} }} = \left| {\frac{{{\mathop{\rm Im}\nolimits} \left[ {{w_{ref}}\left( {x,y} \right)} \right] - {\mathop{\rm Im}\nolimits} \left[ {w\left( {x,y} \right)} \right]}}{{{\mathop{\rm Im}\nolimits} \left[ {{w_{ref}}\left( {x,y} \right)} \right]}}} \right|,
$$
where ${w_{ref}}\left( {x,y} \right)$ is the reference. The highly accurate reference values can be generated by using, for example, Algorithm 680 \cite{Poppe1990a, Poppe1990b}, Algorithm 916 \cite{Zaghloul2011} or a recent C++ code in {\it{RooFit}} package from the CERN's library \cite{Karbach2014}.

Consider Figs. 2 and 3 illustrating logarithms of relative errors of the algorithm for the real and imaginary parts over the domain $0 \le x \le 10 \cap 0 \le y \le 10$, respectively. The rough surface of these plots indicates that the computation reaches the limit of double precision $\sim {10^{ - 16}}$. In particular, while $y > 0.1$ the accuracy can exceed ${10^{ - 15}}$. However, there is a narrow domain along the along $x$ axis in Fig. 2 (red color area) where the accuracy deteriorates by about an order of the magnitude. Apart from this, we can also see in Fig. 3 a sharp peak located at the origin where accuracy is $\sim {10^{ - 14}}$. Figures 4 and 5 depict these areas of the plots magnified for the real and imaginary parts, respectively. 

It should be noted, however, that the worst accuracy $\sim {10^{ - 14}}$ is relatively close to the limitation at double precision computation. Furthermore, since the corresponding areas are negligibly small as compared to the entire inner circle area, where the approximations \eqref{eq_10} and \eqref{eq_14} are applied, their contribution is ignorable and, therefore, does not affect the average accuracy $\sim {10^{ - 15}}$. Thus, the computational test reveals that the obtained accuracy at double precision computation for the complex error function is absolutely consistent with {\it{CERNLIB}}, {\it{libcerf}} and {\it{RooFit}} packages (see the work \cite{Karbach2014} for specific details regarding accuracy of these packages). The run-time test has been compared with the Ab-Initio group of MIT implementation in C/C++ \cite{Johnson, Matlab38787}.

\section{Chiarella and Reichell approximation}

A complex domain with vanishing imaginary part $y=\mathop{\rm Im}\left[z\right]\to 0$ is commonly considered difficult for high accuracy computation of the complex error function. The Chiarella and Reichel approximation (see Appendix D)
\[
\begin{aligned}
w\left( z \right)\approx &\,i\frac{h}{\pi \,z}-i\frac{2h\,z}{\pi }\sum\limits_{n=1}^{\infty}{\frac{{{e}^{-{{n}^{2}}{{h}^{2}}}}}{{{n}^{2}}{{h}^{2}}-{{z}^{2}}}} \\ 
& +\frac{2{{e}^{-{{z}^{2}}}}}{1-{{e}^{-2\pi iz/h}}}\left\{ 1-H\left( y-\frac{\pi }{h} \right) \right\},  
\end{aligned}
\]
where $h>0$ is a small fitting parameter and $H\left(y-\pi/h\right)$ is the Heaviside step function right-shifted by $\pi/h$, may be used to resolve such a problem. However, as we can see the Chiarella and Reichel approximation contains poles at $z^2 = n^2h^2$ located along $x$-axis. Consequently, its application leads to computational error/overflow near vicinities of each poles. Another issue that has to be taken into account is the presence of the Heaviside step function. In particular, the Heaviside step function $H\left( y-\pi/h\right)$ cannot be always taken equal to unity over the entire HITRAN domain $0 \le x \le 40,000 \cap {10^{ - 4}} \le y \le {10^2}$ since such a simplification results in a rapid accuracy deterioration at $y\lesssim 0.5$. This signifies a requirement for additional logical condition and corresponding function file associated with the last term
$$
\frac{2{{e}^{-{{z}^{2}}}}}{1-{{e}^{-2\pi iz/h}}}\left\{ 1-H\left( y-\frac{\pi }{h} \right) \right\}
$$
that, consequently, leads to deceleration of the Matlab code.

\section{Salzer's approximation}

Another approach for computation of the complex error function that covers the complex domain with vanishing imaginary part $y=\mathop{\rm Im}\left[z\right]\to 0$ providing high accuracy is based on the Salzer's approximation. In 1951 Salzer showed that by using the Poisson summation formula
$$
\sum\limits_{n=-\infty }^{\infty }{{{e}^{-{{\left(u+na \right)}^{2}}}}=\frac{\sqrt{\pi }}{a}}\sum\limits_{n=-\infty }^{\infty }{{{e}^{-\frac{{{n}^{2}}{{\pi }^{2}}}{{{a}^{2}}}}}\cos \left( \frac{2n\pi }{a}u \right)},
$$
where $a>0$ is a fitting parameter, the error function of the complex argument can be approximated with high accuracy by using the following approximation \cite{Salzer1951}
\footnotesize
\[
\begin{aligned}
\text{erf}\left( x,y \right)\approx &\,\text{erf}\left( x,0 \right)+\frac{a}{\pi }{{e}^{-{{x}^{2}}}} \\
&\times \left[ f\left( x,y \right)+2\sum\limits_{n=1}^{N}{\frac{{{e}^{-{{a}^{2}}{{n}^{2}}}}}{{{a}^{2}}{{n}^{2}}+{{x}^{2}}}}\left( x-\frac{x\cosh \left( 2any \right)-ian\sinh \left( 2any \right)}{{{e}^{2ixy}}} \right) \right],
\end{aligned}
\]
\normalsize
where $\text{erf}\left( x,0 \right)$ is the error function of the real argument $x$, $a$ is a small parameter and
\small
$$
f\left( x,y \right)=\left\{ 
\begin{aligned}
& \frac{1-{{e}^{-2ixy}}}{x}, &\qquad x\ne 0 \\ 
& 2iy, &\qquad x = 0. \\ 
\end{aligned} 
\right.
$$
\normalsize
This equation is widely used at $a=1/2$ (see for example equation (5) in \cite{Zaker1969}). Substituting this equation into the identity \cite{Abramowitz1972, Schreier1992}
$$
w\left( z \right)={{e}^{-{{z}^{2}}}}\text{erfc}\left( -iz \right)={{e}^{-{{z}^{2}}}}\left( 1-\text{erf}\left( -iz \right) \right)
$$
yields the reformulated approximation for the complex error function
\footnotesize
\begin{equation}
\begin{aligned}\label{eq_15}
w\left( x,y \right)\approx &\,{{e}^{-{{\left( x+iy \right)}^{2}}}}\Bigg\{ 1-\text{erf}\left( y,0 \right)-\frac{a}{\pi }{{e}^{-{{y}^{2}}}}\Bigg[ f\left( y,-x \right) \\
&+\left.\left.2\sum\limits_{n=1}^{N}{\frac{{{e}^{-{{a}^{2}}{{n}^{2}}}}}{{{a}^{2}}{{n}^{2}}+{{y}^{2}}}}\left( y-{{e}^{2ixy}}\left( y\cosh \left( 2anx \right)+ian\sinh \left( 2anx \right) \right) \right) \right] \right\}.
\end{aligned}
\end{equation}
\normalsize

Zaghloul and Ali showed that the real and imaginary parts of the equation \eqref{eq_15} are given by \cite{Zaghloul2011}
\footnotesize
\[
\tag{16a}\label{eq_16a}
\begin{aligned}
\operatorname{Re}\left[ w\left( x,y \right) \right] \approx &\,{{e}^{-{{x}^{2}}}}\Bigg\{ {{e}^{{{y}^{2}}}}\text{erfc}\left( y,0 \right)\cos \left( 2xy \right)+\frac{2a}{\pi } \\
&\times \left.\left[ x\sin \left( xy \right)\text{sinc}\left( xy \right)+y\sum\limits_{n=1}^{N}{\frac{{{e}^{-{{a}^{2}}{{n}^{2}}}}}{{{a}^{2}}{{n}^{2}}+{{y}^{2}}}}\left( \cosh \left( 2anx \right)-\cos \left( 2xy \right) \right) \right] \right\}
\end{aligned}
\]
\normalsize
and
\footnotesize
\[
\tag{16b}\label{eq_16b}
\begin{aligned}
\operatorname{Im}\left[ w\left( x,y \right) \right]\approx & \,{{e}^{-{{x}^{2}}}}\Bigg\{ -{{e}^{{{y}^{2}}}}\text{erfc}\left( y,0 \right)\sin \left( 2xy \right)+\frac{2a}{\pi } \\
&\times\left.\left[ x\,\text{sinc}\left( 2xy \right)+\sum\limits_{n=1}^{N}{\frac{{{e}^{-{{a}^{2}}{{n}^{2}}}}}{{{a}^{2}}{{n}^{2}}+{{y}^{2}}}}\left( y\sin \left( 2xy \right)+an\sinh \left( 2anx \right) \right) \right] \right\},
\end{aligned}
\]
\normalsize
respectively, where the sinc function is defined such that
$$
\left\{ \text{sinc}\left( xy\ne 0 \right)=\sin \left( xy \right)/\left( xy \right),\text{sinc}\left( xy=0 \right)=1 \right\}.
$$
Zaghloul and Ali also noticed the significance of the small parameter $a$ and its dependance on the number of the summation terms $N$. In particular, when we take $a=0.5$ and $N=23$, the worst accuracy of the Salzer's approximation \eqref{eq_15} in the considered domain $0 \le y < 0.1 \cap \left| {x + iy} \right| \le 8$ is $\sim 10^{-8}$. However, at $a=0.56$ and $N=23$ the worst accuracy becomes $\sim 10^{-14}$. This can be seen from the Figs. 7 and 8 showing the real \eqref{eq_16a} and imaginary \eqref{eq_16b} parts of the Salzer's approximation \eqref{eq_15}.

Although the Salzer's approximation \eqref{eq_15} also covers the required domain $0 \le y < 0.1 \cap \left| {x + iy} \right| \le 8$ with high accuracy $\sim 10^{-14}$, the application of the proposed approximation \eqref{eq_14} is more advantageous since the Salzer's approximation \eqref{eq_15} is not a rational function. Thus, according to Berk and Hawes in applications requiring large-size input arrays it is necessary to avoid trigonometric dependencies of the variable $x$ in order to provide more rapid computation \cite{Berk2017}. As we can see, the equations \eqref{eq_15}, \eqref{eq_16b} and \eqref{eq_16b} contain trigonometric functions $\sinh \left( 2anx \right)$ and $\cosh \left( 2anx \right)$ dependent upon the input parameter $x$ and index $n$ running from $1$ to $N$. As a result, this requires $N$ computations of these functions that strongly decelerate the computation.

Since the proposed approximation \eqref{eq_14} combines a simple rational function $2ih{e^{{\sigma ^2}}}z\theta \left( {{z^2} + {\sigma ^2}} \right)$ with only a single non-rational (exponential) term $e^{-z^2} = e^{-(x+iy)^2}$ dependent on the input parameter $x$, it is significantly faster in computation than the Salzer's approximation \eqref{eq_15}. In particular, in approximation \eqref{eq_14} the term $e^{-(x+iy)^2}$ is computed only once, whereas in the Salzer's approximation \eqref{eq_15} the hyperbolic functions $\sinh(2anx)$ and $\cosh(2anx)$ has to be computed $N$ times each while index $n$ runs from $1$ to $N$. Consequently, in terms of computational speed the considered approximation \eqref{eq_14} may be more efficient, for example, in radiative transfer applications like MODTRAN \cite{Berk2017} or bytran \cite{Pliutau2017}, where approximations from our earlier works \cite{Abrarov2011, Abrarov2015f} are currently implemented.

\section{Run-time test}

The Matlab is an array programing language. Consequently, in contrast to the scalar programing languages such as C/C++ or Fortran its run-time is strongly dependent on the logical conditions required for computational flow. Particularly, in a Matlab implementation it is very desirable to reduce the number of the applied approximations in order to minimize the computational flow of the large-size arrays for more rapid calculation. Therefore, in the Matlab code shown in Appendix E we applied three approximations only. These approximations are the Laplace continued fraction \cite{Gautschi1970, Jones1988, Poppe1990a, Poppe1990b}, equation \eqref{eq_10} and the new equation \eqref{eq_14}.

The run-time test has been compared with current version of the Ab-Initio group of MIT implementation developed by Steven Johnson in C/C++ programming language \cite{Johnson, Matlab38787}. This implementation utilizes a large number of approximations and, therefore, cannot be translated directly into Matlab code for rapid computation. However, as a C/C++ implementation it is well-optimized for the rapid and highly accurate  computation of the family of the Faddeeva functions with relative error smaller than $10^{-13}$ \cite{Johnson}. In order to resolve such a problem, Steven Johnson proposed several useful plugins (MEX files) to run C/C++ source files from the Matlab environment \cite{Johnson, Matlab38787}.

Although this technique enables us to integrate C/C++ source files into the Matlab, its practical application is encountered with some complexities such as the requirement for programing experience in C++ to run, script and modify the source files according to the users' needs as well as proper installations of relevant software. Furthermore, there may also be some incompatibility issues between different versions of the Matab distributive and C/C++ compilers (see forum in \cite{Matlab38787} for the possible incompatibilities). This requires renovation of the previously installed Microsoft C++ libraries and applications that, as a consequence, may cause unexpected side effects in running other computer programs. Therefore, our practical experience with researchers specializing in applications of the Voigt/complex error function suggests that the direct Matlab implementation alone rather than its integration with C/C++ source files would be more convenient for majority of the Matlab users. However, it should be noted that the users who are familiar with both C/C++ and Matlab may find such an integration very useful and efficient.

We compared computational speed over a domain that is common for equations \eqref{eq_16a}, \eqref{eq_16b} associated with C/C++ source files {\it{Faddeeva.hh}} and {\it{Faddeeva.cc}} with equation \eqref{eq_14} applied in the Matlab code shown in Appendix E. Thus, for a random array $z=x+iy$ consisting of $10$ million elements in the ranges $0<x<6$ and $0<y<0.1$ at smallest relative error equal to the DBL\_EPSILON, the computational test reveals about $11.6$ and $7.5$ seconds for the C/C++ and Matlab implementations, respectively. For example, the following C/C++ code can be used to verify the computational speed of the Salzer's equations \eqref{eq_16a} and \eqref{eq_16b}:
\footnotesize
\begin{verbatim}
/**
These command lines can be added at the end of the source
file 'Faddeeva.cc' to verify the computational speed of the
Salzer's equations (16a) and (16b).
*/

#ifdef __cplusplus
# include <cstdio>
# include <cstdlib>
# include <iostream>
# include <iomanip>

#else
# include <stdio.h>
#endif

/**
The function fRand(minNum, maxNum) generates a random number
within the range from 'minNum' to 'maxNum'.
*/

double fRand(double minNum, double maxNum){

    double fMin = minNum, fMax = maxNum;
    double f = (double)rand() / RAND_MAX;

return fMin + f * (fMax - fMin);
};

int main(void){

    for(int k = 0; k < 1e7; k++){ // execute the computation ...
    // 10 million times with random numbers 0 < x < 6 and 0 < y < 0.1
    
        FADDEEVA(w)(C(fRand(0,6),fRand(0,0.1)),DBL_EPSILON);
    };
// cout << setprecision(16) << FADDEEVA(w)(C(5.0,1e-6),1e-14) << endl;

return 0;
};
\end{verbatim}
\normalsize

The corresponding Matlab code that verifies the computational speed of the equation \eqref{eq_14} is given by:
\small
\begin{verbatim}
% ******************************************************************
x = 6*rand(1e7,1); % this generates 10 million random numbers ...
                   % within the range 0 < x < 6
y = 0.1*rand(1e7,1); % this generates 10 million random numbers ...
                     % within the range 0 < y < 0.1
z = x + 1i*y; % input array
tic;fadfunc(z);toc % show the run-time
% ******************************************************************
\end{verbatim}
\normalsize

We also compared computational speed over a domain for a random array consisting of $10$ million elements in the ranges $0<x<10,000$ and $0<y<10,000$. For this case, at smallest relative error equal to the DBL\_EPSILON, the computational test reveals about $2.1$ and $2.8$ seconds for the C/C++ and Matlab implementations, respectively.

As we can see, the overall C/C++ performance is more rapid than the Matlab by a factor $1.33$. However, over the specific domain $0 < x < 6$ and $0 < y < 0.1$ that involves the Salzer's equations \eqref{eq_16a}, \eqref{eq_16b} and the proposed equation \eqref{eq_14}, the provided Matlab code is faster by a factor $1.54$. Despite that the Matlab is usually slower than C/C++, this run-time test demonstrates that the new equation \eqref{eq_14} is more rapid in computation than the Salzer's equations \eqref{eq_16a} and \eqref{eq_16b}. Furthermore, the most important ranges for the variables $x$ and $y$ are located within the complex domain $\left|x+iy\right| \lesssim 15$ since, mathematically, the sharp (central) regions of the spectral lines rather than their flat regions are more difficult to interpolate in spectral line fitting. Consequently, there is no any specific reason to distribute equidistantly a grid-point spacing for computation; the regions with larger $x$ or $y$ can be computed with relatively sparse grid-points (and even approximated more rapidly with a simple polynomial, if necessary) since a high spectral resolution for the flat curve regions is not required. It is very common a non-uniform distribution of the grid-points in radiative transfer applications to perform more rapid and efficient computation (see for example \cite{Quine2002}). Therefore, computation at smaller values of the parameters $x$ and $y$ is primarily important practically in all spectroscopic applications. Thus, we can see that more rapid performance (by about $30\%$) of the C/C++ implementation that occurs at larger parameters $x$ and $y$ can be readily overcome when we employ this type of non-uniform distribution of the grid-points because in this case the computational speed of the Matlab code becomes about same or even faster. For example, by replacing command lines in the C/C++ and Matlab codes above as
\small
\begin{verbatim}
    FADDEEVA(w)(C(fRand(0,15),fRand(0,15)),DBL_EPSILON);
\end{verbatim}
\normalsize
and
\small
\begin{verbatim}
    x = 15*rand(1e7,1); % this generates 10 million random numbers ...                   
                        % within the range 0 < x < 15
    y = 15*rand(1e7,1); % this generates 10 million random numbers ...
                        % within the range 0 < y < 15
    z = x + 1i*y; % input array
\end{verbatim}
\normalsize
we obtain in average $5.9$ and $4.0$ seconds, respectively. This signifies that if we take, say, $9$ million grid-points over complex domain $\left|x+iy\right| < 15$ corresponding to the sharp (central) regions of the spectral lines and $1$ million grid-points over complex domain $\left|x+iy\right| \geq 15$ corresponding to the their flat regions, then the Matlab will remain faster in computation than C/C++ implementation.

To the best of our knowledge, the proposed Matlab code is the most rapid in computation as compared to any other Matlab codes with comparable accuracy ever reported in scientific literature or online elsewhere (see, for example, Matlab codes in \cite{Zaghloul2011, Matlab12091, Matlab47801} for computation of the complex error function covering the entire complex domain with high accuracy). This has been archived due to application of only three rapid equations in the Matlab code.

C/C++ source files were compiled by using the Intel C/C++ compiler. The run-time test has been performed on a typical desktop computer Intel(R), CPU at $2.20$ GHz, RAM $8$ GB on the Windows $10$ platform.

\section{Conclusion}

A rational approximation of the Dawson\text{'s} integral is derived and implemented for rapid and highly accurate computation of the complex error function. The computational test we performed shows the accuracy exceeding $\sim {10^{ - 14}}$ in the domain of practical importance $0 \le y < 0.1 \cap \left| {x + iy} \right| \le 8$. A Matlab code for computation of the complex error function covering the entire complex plane is presented.

\section*{Acknowledgments}
This work is supported by National Research Council Canada, Thoth Technology Inc. and York University. The authors wish to thank the principal developer of the MODTRAN Dr. Alexander Berk for constructive discussions and useful suggestions.

\section*{Appendix A}

It is not difficult to show that the complex error function \eqref{eq_1} can be expressed alternatively as (see equation (3) in \cite{Srivastava1987} and \cite{Srivastava1992}, see also Appendix A in \cite{Abrarov2015a} for derivation)
\[
\tag{A.1}\label{eq_A.1}
w\left( {x,y} \right) = \frac{1}{{\sqrt \pi  }}\int\limits_{ - \infty }^\infty  {\exp \left( { - {t^2}/4} \right)\exp \left( { - yt} \right)\exp \left( {ixt} \right)dt}.
\]
Consequently, from this equation we have
\[\tag{A.2}\label{eq_A.2}
\begin{aligned}
w\left( {x,y} \right) + w\left( { - x, - y} \right) &= \frac{1}{{\sqrt \pi  }}\int\limits_{ - \infty }^\infty  {{e^{ - {t^2}/4}}\left[ {{e^{\left( { - y + ix} \right)t}} + {e^{\left( {y - ix} \right)t}}} \right]dt} \\
 &= \frac{1}{{\sqrt \pi  }}\int\limits_{ - \infty }^\infty  {{e^{ - {t^2}/4}}\left[ {{e^{i\left( {x + iy} \right)t}} + {e^{ - i\left( {x + iy} \right)t}}} \right]dt}.
\end{aligned}
\]
Using the Euler\text{'}s identity
$$
\frac{{{e^{i\mu }} + {e^{ - i\mu }}}}{2} = \cos \left(\mu\right),
$$
where $\mu \in \mathbb{C}$, from the equation \eqref{eq_A.2} it follows that
$$
w\left( {x,y} \right) + w\left( { - x, - y} \right) = \frac{2}{{\sqrt \pi  }}\int\limits_{ - \infty }^\infty  {{e^{ - {t^2}/4}}\cos \left( {\left[ {x + iy} \right]t} \right)dt}  = 2{e^{ - {{\left( {x + iy} \right)}^2}}}
$$
or
$$
w\left( z \right) + w\left( { - z} \right) = 2{e^{ - {z^2}}}.
$$

\section*{Appendix B}

The forward Fourier transform can be defined by following equation \cite{Bracewell2000, Hansen2014}
\[
\label{eq_B.1}\tag{B.1}
F\left( \nu  \right) = {\cal{F}}\left\{ {f\left( t \right)} \right\}\left( \nu  \right) = \int\limits_{ - \infty }^\infty  {f\left( t \right){e^{ - 2\pi i \nu t}}dt}.
\]
It has been shown in our recent publication that sampling with the Gaussian function of the kind $h{e^{ - {{\left( {t/c} \right)}^2}}}/ \left( {c\sqrt \pi  } \right)$, where $h$ is the step between two adjacent sampling points and $c$ is the fitting parameter, leads to the approximation for the inverse Fourier transform as given by \cite{Abrarov2015e}
\[
\label{eq_B.2}\tag{B.2}
f\left( t \right) = {\cal{F}}^{-1}\left\{ {F\left( \nu  \right)} \right\}\left( t \right) \approx h{e^{ - {{\left( {\pi ct} \right)}^2}}}\sum\limits_{n =  - N}^N {F\left( {nh} \right){e^{2\pi itnh}}}.
\]

The equation in \eqref{eq_B.2} is particularly efficient to approximate a non-periodic solitary wavelet (or pulse) due to damping caused by the exponential multiplier ${e^{ - {{\left( {\pi ct} \right)}^2}}}$. Although this damping function excludes periodicity, its presence is not desirable for numerical integration. At $c=0$ the damping effect from the function ${e^{ - {{\left( {\pi ct} \right)}^2}}}$ disappear since it becomes equal to one. As a consequence, the right side of the equation \eqref{eq_B.2} appears to be periodic. Therefore, because of periodicity occurring at $c=0$ this equation becomes valid only within the period $- 1/\left(2h\right) \le t \le 1/\left(2h\right)$. Thus, for this specific case the equation \eqref{eq_B.2} can be rewritten in form of the Fourier series
\[
\label{eq_B.3}\tag{B.3}
f\left( t \right) \approx {h}\sum\limits_{n =  - N}^N {F\left( {n{h}} \right){e^{2\pi itn{h}}}}, \qquad - \frac{1}{{2{h}}} \le t \le \frac{1}{{2{h}}}.
\]

Consider now the function $f\left( t \right) = {e^{ - {t^2}/4}}$. The forward Fourier transform of this function can be found analytically by substituting it into equation \eqref{eq_B.1}. These leads to
$$F\left( \nu  \right) = \int\limits_{ - \infty }^\infty  {{e^{ - {t^2}/4}}{e^{ - 2\pi ivt}}dt}  = 2\sqrt \pi  {e^{ - {{\left( {2\pi \nu } \right)}^2}}}.
$$
Substituting $F\left( \nu \right) = 2\sqrt \pi  {e^{ - {{\left( {2\pi \nu } \right)}^2}}}$ into equation \eqref{eq_B.3} yields the following approximation for the exponential function
\[
{e^{ - {t^2}/4}} \approx 2\sqrt \pi  {h}\sum\limits_{n =  - N}^N {{e^{ - {{\left( {2\pi n{h}} \right)}^2}}}{e^{2\pi itn{h}}}},		\qquad - \frac{1}{{2{h}}} \le t \le \frac{1}{{2{h}}},
\]
or
\[
\label{eq_B.4}\tag{B.4}
{e^{ - {t^2}/4}} \approx 2\sqrt \pi  {h}\sum\limits_{n =  - N}^N {{e^{ - {{\left( {2\pi n{h}} \right)}^2}}}\cos \left( {2\pi tn{h}} \right)},	\qquad - \frac{1}{{2{h}}} \le t \le \frac{1}{{2{h}}}.
\]
Since
$$
{e^{ - {{\left( {2\pi 0{h}} \right)}^2}}}\cos \left( {2\pi t0{h}} \right) = 1
$$
and
$$
\sum\limits_{n =  - N}^{ - 1} {{e^{ - {{\left( {2\pi n{h}} \right)}^2}}}\cos \left( {2\pi tn{h}} \right)}  = \sum\limits_{n = 1}^N {{e^{ - {{\left( {2\pi n{h}} \right)}^2}}}\cos \left( {2\pi tn{h}} \right)}
$$
the approximation \eqref{eq_B.4} can be further simplified as
\small
\[
\label{eq_B.5}\tag{B.5}
{e^{ - {t^2}/4}} \approx 2\sqrt \pi  {h}\left[ {1 + 2\sum\limits_{n = 1}^N {{e^{ - {{\left( {2\pi n{h}} \right)}^2}}}\cos \left( {2\pi tn{h}} \right)} } \right],	\qquad - \frac{1}{{2{h}}} \le t \le \frac{1}{{2{h}}}.
\]
\normalsize

The limitation on the right side $t \le 1/\left( {2{h}} \right)$ along the positive $t$-axis in equation \eqref{eq_B.5} can be excluded by multiplying its left and right sides to $\exp \left( { - \sigma t} \right)$ if the positive constant $\sigma$ is chosen to be sufficiently large. This is possible to achieve since the exponential multiplier $\exp \left( { - \sigma t} \right)$ acts as a dumping function at larger values of the argument $t > 1/\left(2h\right)$. Consequently, if the constant $\sigma $ is large enough, say approximately equal or greater than $1$, then we can write the approximation
$$
{e^{ - {t^2}/4}}{e^{ - \sigma t}} \approx 2\sqrt \pi  {h}\left[ {1 + 2\sum\limits_{n = 1}^N {{e^{ - {{\left( {2\pi n{h}} \right)}^2}}}\cos \left( {2\pi tn{h}} \right)} } \right]{e^{ - \sigma t}}, \qquad \sigma  \mathbin{\lower.3ex\hbox{$\buildrel>\over
{\smash{\scriptstyle\sim}\vphantom{_x}}$}} 1,
$$
that remains always valid without any limitation along the positive $t$-axis. Therefore, taking $y \ge 0$ we get
\footnotesize
\[
\label{eq_B.6}\tag{B.6}
{e^{ - {t^2}/4}}{e^{ - \left( {y + \sigma } \right)t}} \approx 2\sqrt \pi  {h}\left[ {1 + 2\sum\limits_{n = 1}^N {{e^{ - {{\left( {2\pi n{h}} \right)}^2}}}\cos \left( {2\pi tn{h}} \right)} } \right]{e^{ - \left( {y + \sigma } \right)t}}, \qquad \sigma  \mathbin{\lower.3ex\hbox{$\buildrel>\over
{\smash{\scriptstyle\sim}\vphantom{_x}}$}} 1.
\]
\normalsize
Since ${e^{ - {t^2}/4}}{e^{ - yt}} = {e^{{\sigma ^2}}}{e^{ - {{\left( {t - 2\sigma } \right)}^2}/4}}{e^{ - \left( {y + \sigma } \right)t}}$ from approximation \eqref{eq_B.6} we obtain
\footnotesize
\[
\label{eq_B.7}\tag{B.7}
{e^{ - {t^2}/4}}{e^{ - yt}} \approx 2\sqrt \pi  {h}{e^{{\sigma ^2}}}\left[ {1 + 2\sum\limits_{n = 1}^N {{e^{ - {{\left( {2\pi n{h}} \right)}^2}}}\cos \left( {2\pi n{h}\left( {t - 2\sigma } \right)} \right)} } \right]{e^{ - \left( {y + \sigma } \right)t}}, \quad \sigma  \mathbin{\lower.3ex\hbox{$\buildrel>\over
{\smash{\scriptstyle\sim}\vphantom{_x}}$}} 1.
\]
\normalsize
Lastly, substituting approximation \eqref{eq_B.7} into integral \eqref{eq_A.1} we obtain the approximation \eqref{eq_9} for the complex error function.

\section*{Appendix C}

As we have shown in our publication \cite{Abrarov2015d}, the real part of the complex error function can be found as
\[\tag{C.1}\label{eq_C.1}
{\mathop{\rm Re}\nolimits} \left[ {w\left( {x,y} \right)} \right] = \frac{{w\left( {x,y} \right) + w\left( { - x,y} \right)}}{2}
\]
and since \cite{Armstrong1967}
$$
{\mathop{\rm Re}\nolimits} \left[ {w\left( {x,y = 0} \right)} \right] = K\left( {x,y = 0} \right) = {e^{ - {x^2}}},
$$
the substitution of the approximation \eqref{eq_9} at $y = 0$ into equation \eqref{eq_C.1} leads to
\[\tag{C.2}\label{eq_C.2}
\begin{aligned}
{e^{ - {x^2}}} \approx &\frac{{2{e^{{\sigma ^2}}}h\sigma }}{{{x^2} + {\sigma ^2}}} \\
&+ \frac{1}{2}\sum\limits_{n = 1}^N {\left( {\frac{{{A_n} - i\left( {x + i\sigma } \right){B_n}}}{{C_n^2 - {{\left( {x + i\sigma } \right)}^2}}} + \frac{{{A_n} - i\left( { - x + i\sigma } \right){B_n}}}{{C_n^2 - {{\left( { - x + i\sigma } \right)}^2}}}} \right)}.
\end{aligned}
\]

According to the definition of the $L$-function \eqref{eq_7}, at $y = 0$ we can write
\[\tag{C.3}\label{eq_C.3}
L\left( {x,y = 0} \right) = \mathop {\lim }\limits_{y \to 0} \frac{1}{\pi }\int\limits_{ - \infty }^\infty  {\frac{{{e^{ - {t^2}}}\left( {x - t} \right)}}{{{y^2} + {{\left( {x - t} \right)}^2}}}dt}.
\]
Consequently, substituting the approximation for the exponential function \eqref{eq_B.2} into equation \eqref{eq_C.3} yields
\small
\[\tag{C.4}\label{eq_C.4}
\begin{aligned}
L\left( {x,y = 0} \right) \approx &\mathop {\lim }\limits_{y \to 0} \frac{1}{\pi }\int\limits_{ - \infty }^\infty  \frac{{x - t}}{{{y^2} + {{\left( {x - t} \right)}^2}}} \left[ \frac{{2{e^{{\sigma ^2}}}h\sigma }}{{{t^2} + {\sigma ^2}}} \right. \\
&+ \left. \frac{1}{2}\sum\limits_{n = 1}^N {\left( {\frac{{{A_n} - i\left( {t + i\sigma } \right){B_n}}}{{C_n^2 - {{\left( {t + i\sigma } \right)}^2}}} + \frac{{{A_n} - i\left( { - t + i\sigma } \right){B_n}}}{{C_n^2 - {{\left( { - t + i\sigma } \right)}^2}}}} \right)}  \right]dt.
\end{aligned}
\]
\normalsize
The integrand of the integral \eqref{eq_C.4} is analytic everywhere except $2 + 2N$ isolated points on the upper half plane
$$
{t_r} = \left\{ {x + iy,i\sigma , - {C_n} + i\sigma ,{C_n} + i\sigma } \right\}, \qquad n \in \left\{ {1,2,3,\,\, \ldots \,\,N} \right\}.
$$
Therefore, using the Residue Theorem\text{'}s formula
$$
\frac{1}{{2\pi i}}\oint\limits_{{C_{ccw}}} {f\left( t \right)dt = \sum\limits_{r = 1}^{2 + 2N} {{\rm{Res}}\left[ {f,{t_r}} \right]} },
$$
where ${C_{ccw}}$ denotes a contour in counterclockwise direction enclosing the upper half plain (for example as a semicircle with infinite radius) and
\begin{center}
$f\left( t \right) = \frac{1}{\pi } \, \times$ integrand of the integral \eqref{eq_C.4},
\end{center}
we obtain the rational approximation \eqref{eq_11} of the $L$-function.

\section*{Appendix D}

In 1968 Chiarella and Reichel found an elegant series expansion for the function
$$
{{W}_{0}}\left( x,t \right)={{U}_{0}}\left( x,t \right)+i{{V}_{0}}\left( x,t \right)=\frac{\Omega\left(x,t\right)}{{{\left( 4\pi t \right)}^{1/2}}}\int\limits_{-\infty }^{\infty }{\frac{{{e}^{-{{u}^{2}}}}}{{{u}^{2}}+{{\Omega^{2}\left(x,t\right)}}}d}u
$$
that can be represented in form \cite{Chiarella1968, Matta1971}
\[
\tag{D.1}\label{D.1}
\begin{aligned}
 {{W}_{0}}\left( x,t \right)= & \frac{h}{\Omega\left(x,t\right){{\left( 4\pi t \right)}^{1/2}}}+\frac{2h\Omega\left(x,t\right)}{{{\left( 4\pi t \right)}^{1/2}}}\sum\limits_{n=1}^{\infty }{\frac{{{e}^{-{{n}^{2}}{{h}^{2}}}}}{{{\Omega^{2}\left(x,t\right)}}+{{n}^{2}}{{h}^{2}}}} \\ 
 & +\frac{\pi {{e}^{{{\Omega^{2}\left(x,t\right)}}}}}{{{\left( \pi t \right)}^{1/2}}\left( 1-{{e}^{2\pi \Omega\left(x,t\right)/h}} \right)}H\left( t-\frac{{{h}^{2}}}{{{\pi }^{2}}} \right)-\frac{\Omega\left(x,t\right)}{{{\left( 4\pi t \right)}^{1/2}}}E\left( h,\Omega \right),  
\end{aligned}
\]
where ${U}_{0}\left( x,t \right)$ and $V_0\left( x,t \right)$ are the real and imaginary parts of this function, $h>0$ is a small fitting parameter, $\Omega\left(x,t\right) = \left(1-ix\right)/\left(2t^{1/2}\right)$,
$$
H\left( t \right)=\left\{
\begin{aligned}
1, &\quad t > 0 \\
1/2, &\quad t = 0 \\
0, &\quad t < 0
\end{aligned}
\right.
$$
is the Heaviside step function and
\[
E\left( h,\Omega \right)=2{{e}^{-{{\pi }^{2}}/{{h}^{2}}}}\int\limits_{-\infty }^{\infty }{\frac{{{e}^{-{{z}^{2}}}}}{\left\{ {{\left( z-i\pi /h \right)}^{2}}+{{\Omega^{2}\left(x,t\right)}} \right\}\left\{ 1-{{e}^{-2\pi iz/h-2{{\pi }^{2}}/{{h}^{2}}}} \right\}}dz}.
\]

Applying the following identities for the real and imaginary parts
$$
K\left( x,y \right)=\frac{1}{y\sqrt{\pi }}{{U}_{0}}\left( \frac{x}{y},\frac{1}{4{{y}^{2}}} \right)
$$ 
and 
$$
L\left( x,y \right)=\frac{1}{y\sqrt{\pi }}{{V}_{0}}\left( \frac{x}{y},\frac{1}{4{{y}^{2}}} \right)
$$ 
we can reformulate the series expansion \eqref{D.1} in terms of the complex error function according to \eqref{eq_8} as
\[
\tag{D.2}\label{D.2}
\begin{aligned}
w\left( z \right)= &\,i\frac{h}{\pi \,z}-i\frac{2h\,z}{\pi }\sum\limits_{n=1}^{\infty}{\frac{{{e}^{-{{n}^{2}}{{h}^{2}}}}}{{{n}^{2}}{{h}^{2}}-{{z}^{2}}}} \\ 
& +\frac{2{{e}^{-{{z}^{2}}}}}{1-{{e}^{-2\pi iz/h}}}\left\{ 1-H\left( y-\frac{\pi }{h} \right) \right\}+i\frac{z}{\pi }\varepsilon \left( h,z \right),  
\end{aligned}
\]
where
\[
\varepsilon \left( h,z \right)=2{{e}^{-{{\pi }^{2}}/{{h}^{2}}}}\int\limits_{-\infty }^{\infty }{\frac{{{e}^{-{{u}^{2}}}}}{\left\{ {{\left( u-i\pi /h \right)}^{2}}-{{z}^{2}} \right\}\left\{ 1-{{e}^{-2\pi iu/h-2{{\pi }^{2}}/{{h}^{2}}}} \right\}}du}.
\]
Since at $0 < h < 1$ the multiplier ${e}^{-{{\pi }^{2}}/{{h}^{2}}} \approx 0$, the last term $iz\,\varepsilon \left( h,z \right)/\pi$ in equation \eqref{D.2} can be ignored due to its negligibly small contribution.

\bigskip
\section*{Appendix E}
\footnotesize
\begin{verbatim}
function FF = fadfunc(z)

% This program file computes the complex error function, also known as the
% Faddeeva function. It provides high-accuracy and covers the entire
% complex plane. The inner part of the circle |x + i*y| <= 8 is covered by
% the equations (10) and (14). Derivation of the rational approximation
% (10) and its detailed description can be found in the paper [1]. The new
% approximation (14) shown in this paper computes the complex error
% function at small Im[z] << 1. The outer part of the circle |x + i*y| > 8
% is covered by the Laplace continued fraction [2]. The accuracy of this
% function file can be verified by using C/C++ code provided in work [3].
%
% REFERENCES
% [1] S. M. Abrarov and B. M. Quine, A new application methodology of the
%     Fourier transform for rational approximation of the complex error
%     function, arXiv:1511.00774.
%     http://arxiv.org/abs/1511.00774
%
% [2] W. Gautschi, Efficient computation of the complex error function,
%     SIAM J. Numer. Anal., 7 (1970) 187-198.
%     http://www.jstor.org/stable/2949591
%
% [3] T. M. Karbach, G. Raven and M. Schiller, Decay time integrals in
%     neutral meson mixing and their efficient evaluation,
%     arXiv:1407.0748v1 (2014).
%     http://arxiv.org/abs/1407.0748
%
% The code is written by Sanjar M. Abrarov and Brendan M. Quine, York
% University, Canada, December 2015.
%
% Corresponding author's email: absanj AT gmail DOT com

% *************************************************************************
% All parameters in this table are global and can be used anywhere inside
% the program body.
% -------------------------------------------------------------------------
n = 1:23; % define a row vector
sigma = 1.5; % the shift constant
h = 6/(2*pi)/23; % this is the step
% -------------------------------------------------------------------------
% Expansion coefficients for eq. (10)
% -------------------------------------------------------------------------
An = 8*pi*h^2*n.*exp(sigma^2 - (2*pi*h*n).^2).*sin(4*pi*h*n*sigma);
Bn = 4*h*exp(sigma^2 - (2*pi*h*n).^2).*cos(4*pi*h*n*sigma);
Cn = 2*pi*h*n;
% -------------------------------------------------------------------------
% Expansion coefficients for eq. (14)
% -------------------------------------------------------------------------
alpha = 8*pi*h*n*sigma.*exp(-(2*pi*h*n).^2).*sin(4*pi*h*n*sigma);
beta = 2*exp(-(2*pi*h*n).^2).*cos(4*pi*h*n*sigma);
gamma = (2*pi*h*n).^2;
% End of the table
% *************************************************************************

neg = imag(z)<0; % if imag(z) values are negative, then ...
z(neg) = conj(z(neg)); % ... take the complex conjugate
int  = abs(z)<=8; % internal indices

FF = zeros(size(z)); % define array
FF(~int) = extD(z(~int)); % external domain
FF(int) = intD(z(int)); % internal domain

    function ext_d = extD(z) % the Laplace continued fraction

        coeff = 1:11; % define a row vector
        coeff = coeff/2;

        ext_d = coeff(end)./z; % start computing using the last coeff
        for m = 1:length(coeff) - 1
            ext_d = coeff(end-m)./(z - ext_d);
        end
        ext_d = 1i/sqrt(pi)./(z - ext_d);
    end

    function int_d = intD(z) % internal domain

        secD = imag(z)<0.1; % secondary subdomain indices
        int_d(secD) = exp(-z(secD).^2) + 2i*h*exp(sigma^2)*z(secD).* ...
            theta(z(secD).^2 + sigma^2); % compute using eq. (14)
        int_d(~secD) = rAppr(z(~secD)+1i*sigma); % compute using eq. (10)

        function r_appr = rAppr(z) % the rational approximation (10)

            zz = z.^2; % define repeating array

            r_appr = 2i*h*exp(sigma^2)./z;
            for m = 1:23
                r_appr = r_appr + (An(m) - 1i*z*Bn(m))./ ...
                    (Cn(m)^2 - zz);
            end
        end

        function ThF = theta(z) % theta function (see the eq. (14))

            ThF = 1./z;
            for k = 1:23
                ThF = ThF + (alpha(k) + beta(k)*(z - gamma(k)))./ ...
                    (4*sigma^2*gamma(k) + (gamma(k) - z).^2);
            end
        end
    end

%  For negative imag(z) use the identity w(-z) = 2*exp(-z^2) - w(z)
FF(neg) = conj(2*exp(-z(neg).^2) - FF(neg));
end
\end{verbatim}
\normalsize

\bigskip

\newpage
\begin{figure}[ht]
\begin{center}
\includegraphics[width=22pc]{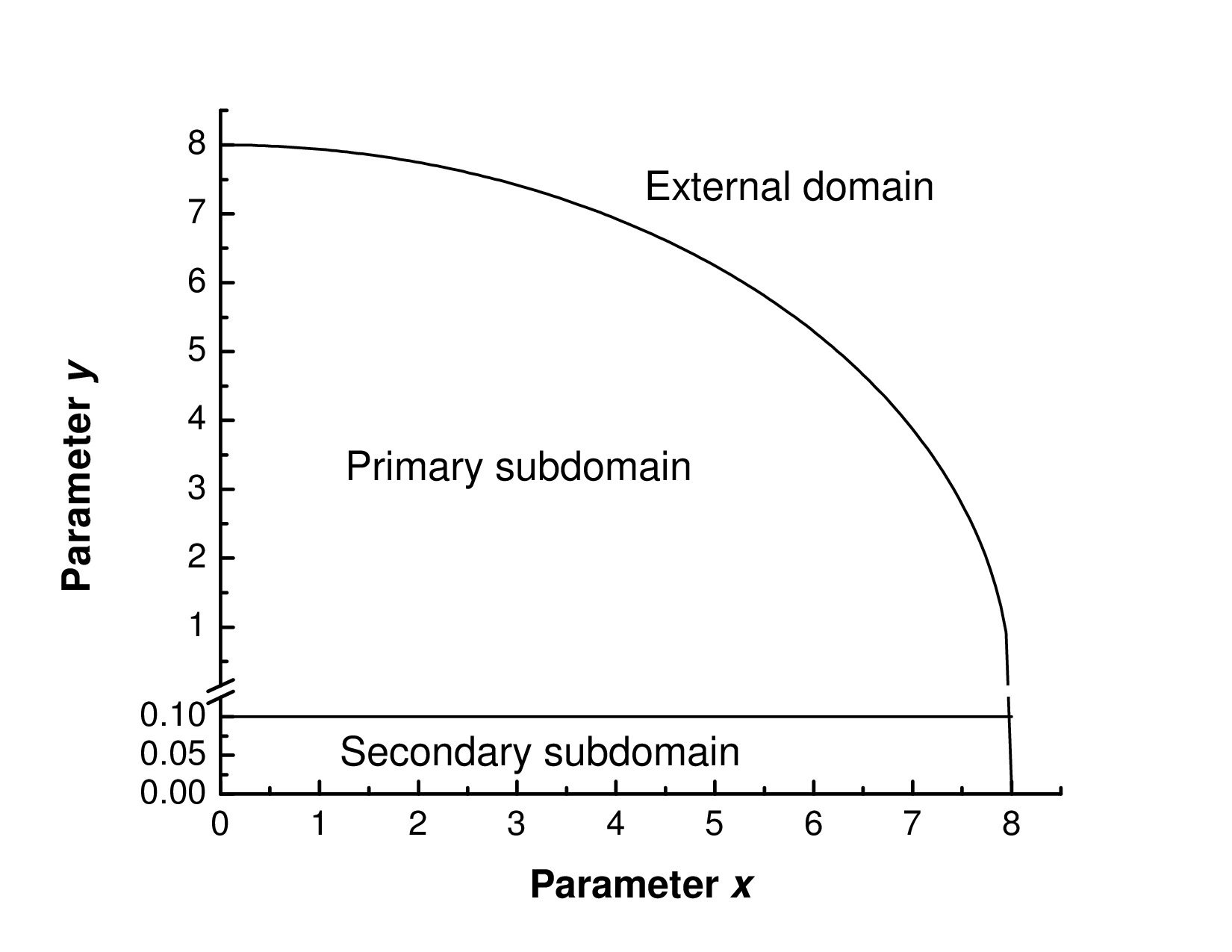}\hspace{1pc}%
\begin{minipage}[b]{28pc}
\vspace{0.3cm}
{\sffamily {\bf{Fig. 1.}} Domains separating the complex plain for high-accuracy computation of the complex error function.}
\end{minipage}
\end{center}
\end{figure}

\newpage
\begin{figure}[ht]
\begin{center}
\includegraphics[width=22pc]{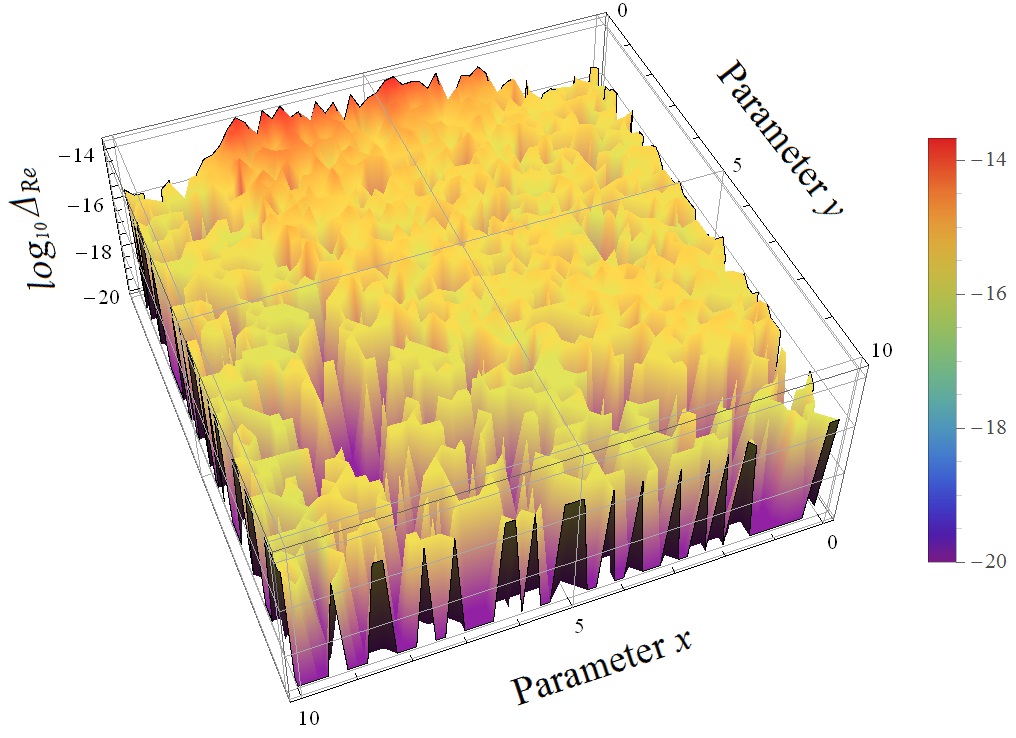}\hspace{1pc}%
\begin{minipage}[b]{28pc}
\vspace{0.3cm}
{\sffamily {\bf{Fig. 2.}} The logarithm of the relative error ${\log _{10}}{\Delta _{{\mathop{\rm Re}\nolimits} }}$ for the real part of the algorithm over the domain $0 \le x \le 10 \cap 0 \le y \le 10$.}
\end{minipage}
\end{center}
\end{figure}

\newpage
\begin{figure}[ht]
\begin{center}
\includegraphics[width=22pc]{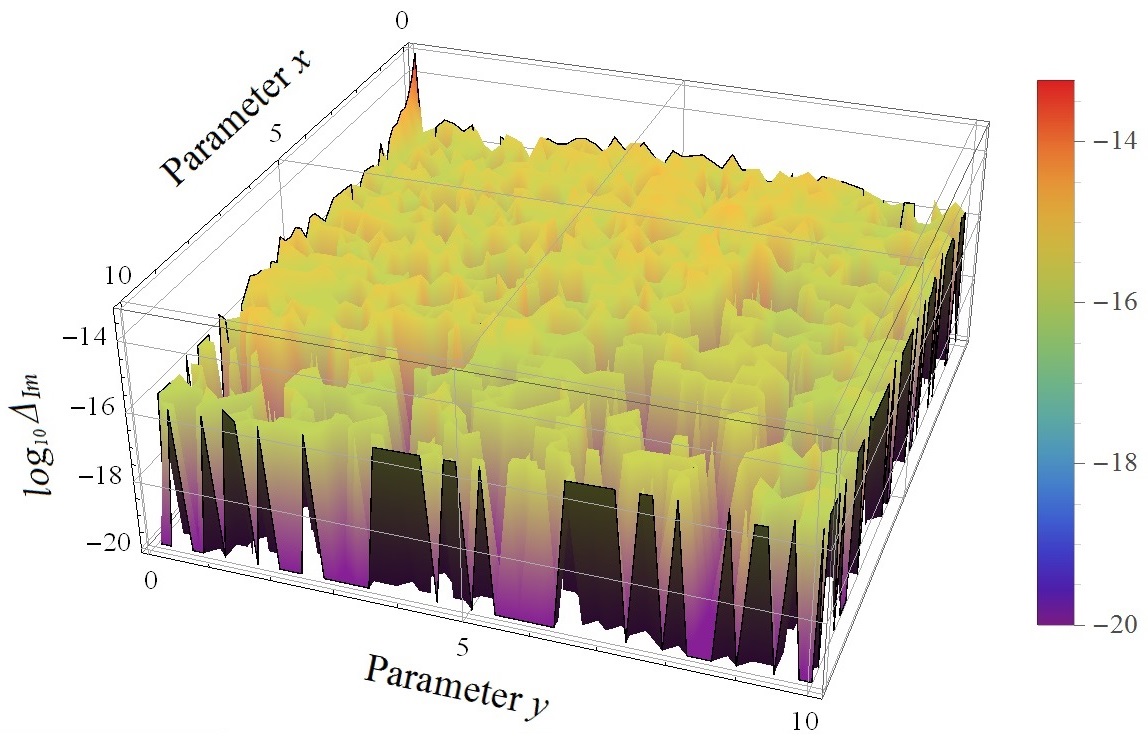}\hspace{1pc}%
\begin{minipage}[b]{28pc}
\vspace{0.3cm}
{\sffamily {\bf{Fig. 3.}} The logarithm of the relative error ${\log _{10}}{\Delta _{{\mathop{\rm Im}\nolimits} }}$ for the imaginary part of the algorithm over the domain $0 \le x \le 10 \cap 0 \le y \le 10$.}
\end{minipage}
\end{center}
\end{figure}

\newpage
\begin{figure}[ht]
\begin{center}
\includegraphics[width=22pc]{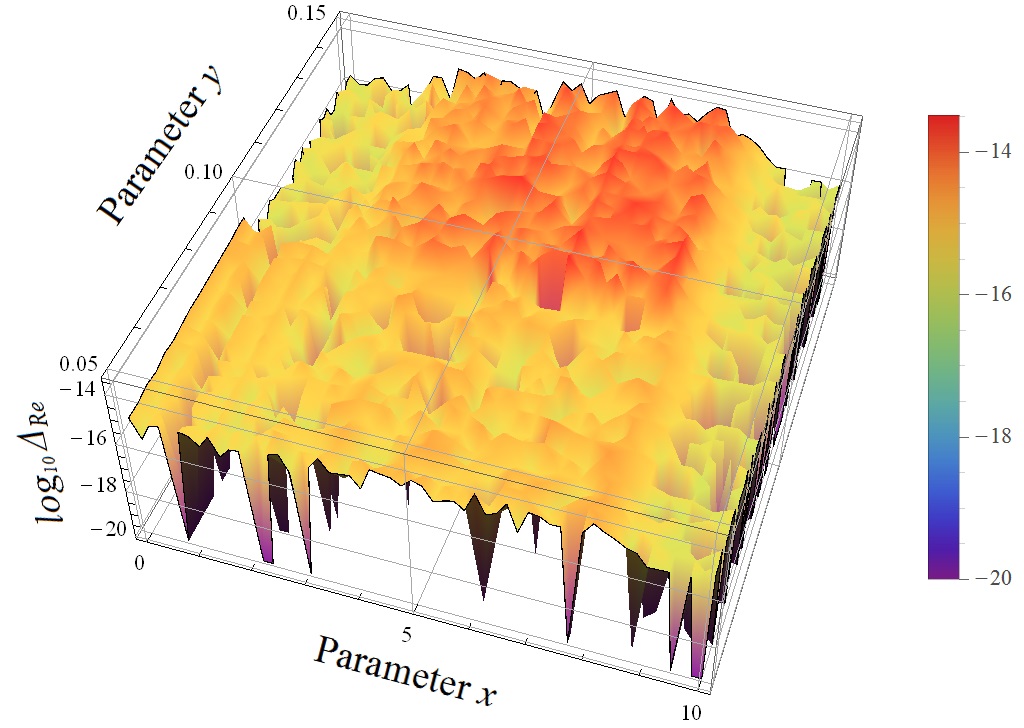}\hspace{1pc}%
\begin{minipage}[b]{28pc}
\vspace{0.3cm}
{\sffamily {\bf{Fig. 4.}} The logarithm of the relative error ${\log _{10}}{\Delta _{{\mathop{\rm Re}\nolimits} }}$ for the real part of the algorithm showing the area with worst accuracy $\sim {10^{ - 14}}$.}
\end{minipage}
\end{center}
\end{figure}

\newpage
\begin{figure}[ht]
\begin{center}
\includegraphics[width=22pc]{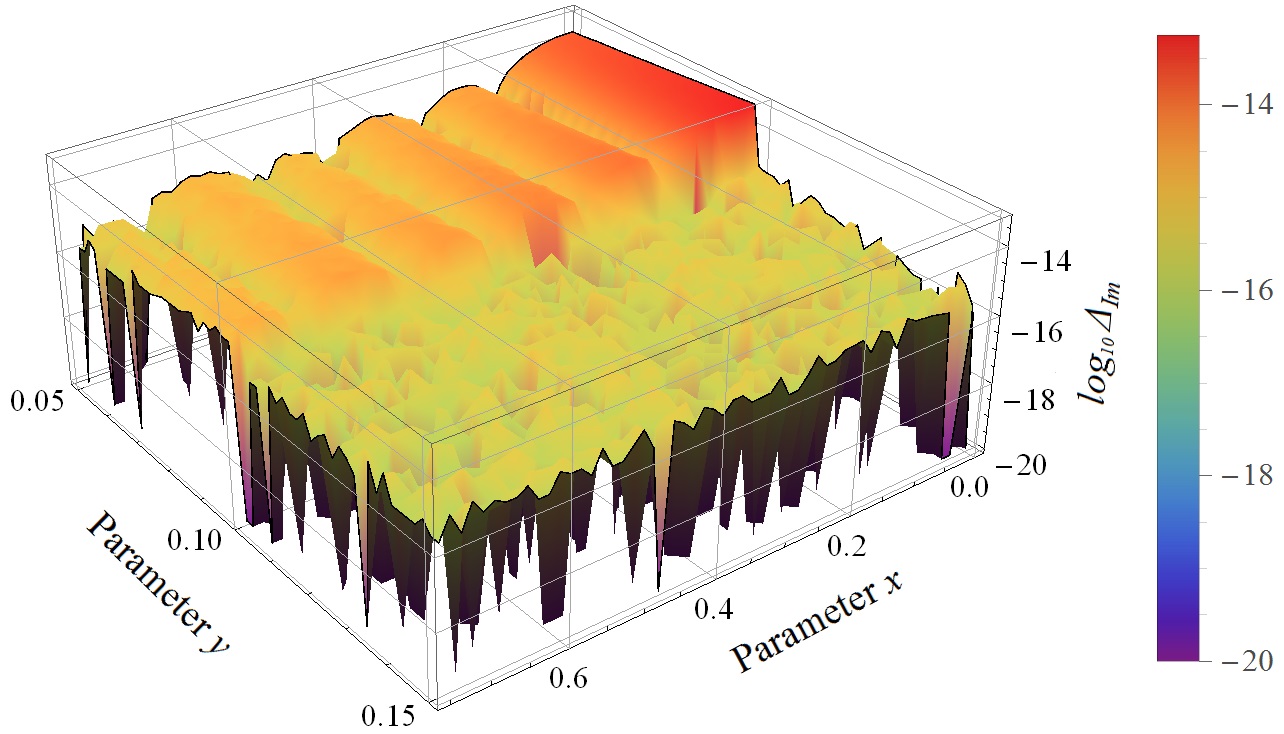}\hspace{1pc}%
\begin{minipage}[b]{28pc}
\vspace{0.3cm}
{\sffamily {\bf{Fig. 5.}} The logarithm of the relative error ${\log _{10}}{\Delta _{{\mathop{\rm Im}\nolimits} }}$ for the imaginary part of the algorithm showing the area with worst accuracy $\sim {10^{ - 14}}$.}
\end{minipage}
\end{center}
\end{figure}

\newpage
\begin{figure}[ht]
\begin{center}
\includegraphics[width=22pc]{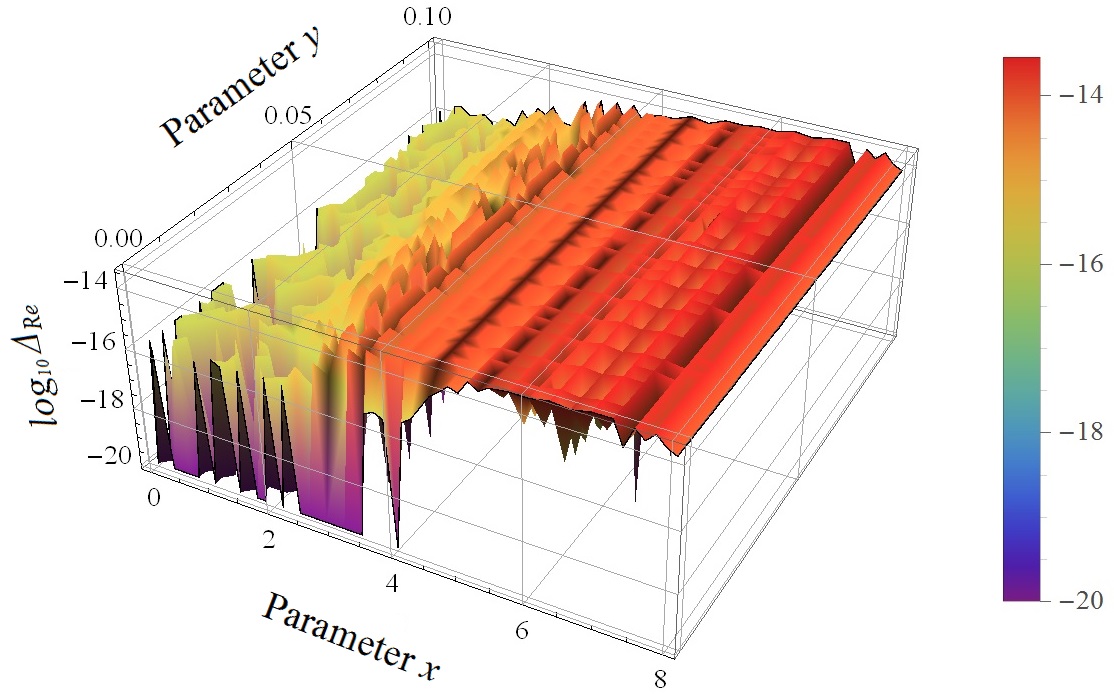}\hspace{1pc}%
\begin{minipage}[b]{28pc}
\vspace{0.3cm}
{\sffamily {\bf{Fig. 6.}} The logarithm of the relative error ${\log _{10}}{\Delta _{{\mathop{\rm Im}\nolimits} }}$ for the equation \eqref{eq_16a} at $a=0.56$ and $N=23$ over the domain $0 \le y < 0.1 \cap \left| {x + iy} \right| \le 8$.}
\end{minipage}
\end{center}
\end{figure}

\newpage
\begin{figure}[ht]
\begin{center}
\includegraphics[width=22pc]{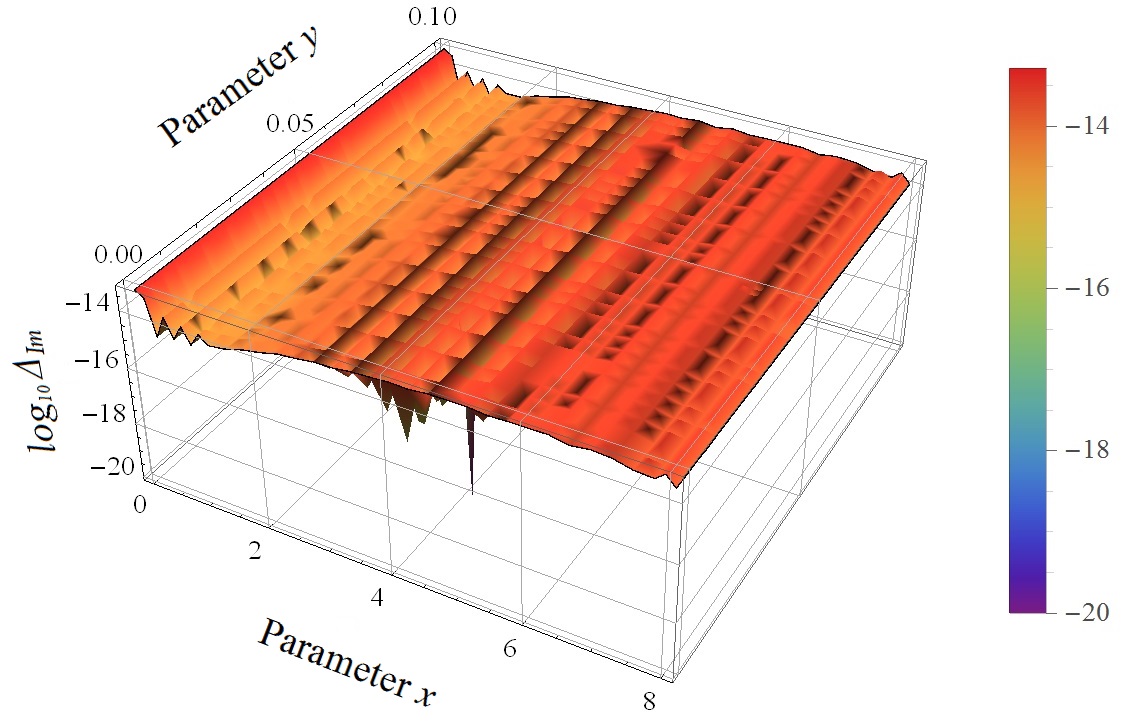}\hspace{1pc}%
\begin{minipage}[b]{28pc}
\vspace{0.3cm}
{\sffamily {\bf{Fig. 7.}}  The logarithm of the relative error ${\log _{10}}{\Delta _{{\mathop{\rm Im}\nolimits} }}$ for the equation \eqref{eq_16b} at $a=0.56$ and $N=23$ over the domain $0 \le y < 0.1 \cap \left| {x + iy} \right| \le 8$.}
\end{minipage}
\end{center}
\end{figure}
\end{document}